\pdfoutput=1
\RequirePackage{ifpdf}
\ifpdf 
\documentclass[pdftex]{sigma}
\else
\documentclass{sigma}
\fi

\newcommand{\bfx}{{\bf x}}
\newcommand{\bfxp}{{{\bf x}^\prime}}
\newcommand{\N}{{\mathbb N}}

\newcommand{\R}{{\mathbb R}}

\newcommand{\Z}{{\mathbb Z}}
\newcommand{\C}{{\mathbb C}}

\begin{document}

\allowdisplaybreaks

\renewcommand{\thefootnote}{$\star$}

\renewcommand{\PaperNumber}{077}

\FirstPageHeading

\ShortArticleName{Def\/inite Integrals using Orthogonality and Integral Transforms}

\ArticleName{Def\/inite Integrals using Orthogonality\\ and Integral Transforms\footnote{This
paper is a contribution to the Special Issue ``Superintegrability, Exact Solvability, and Special Functions''. The full collection is available at \href{http://www.emis.de/journals/SIGMA/SESSF2012.html}{http://www.emis.de/journals/SIGMA/SESSF2012.html}}}

\Author{Howard S.~COHL~$^\dag$ and Hans VOLKMER~$^\ddag$}

\AuthorNameForHeading{H.S.~Cohl and H.~Volkmer}

\Address{$^\dag$~Applied and Computational Mathematics Division,
Information Technology Laboratory,\\
\hphantom{$^\dag$}~National Institute of Standards and Technology,
Gaithersburg, MD, 20899-8910, USA}
\EmailD{\href{mailto:howard.cohl@nist.gov}{howard.cohl@nist.gov}}
\URLaddressD{\url{http://hcohl.sdf.org}}

\Address{$^\ddag$~Department of Mathematical Sciences,
University of Wisconsin-Milwaukee,\\
\hphantom{$^\ddag$}~P.O.~Box 413, Milwaukee, WI, 53201, USA}
\EmailD{\href{mailto:volkmer@uwm.edu}{volkmer@uwm.edu}}

\ArticleDates{Received July 31, 2012, in f\/inal form October 15, 2012; Published online October 19, 2012}

\Abstract{We obtain def\/inite integrals for products of associated Legendre functions
with Bessel functions, associated Legendre functions, and Chebyshev polynomials
of the f\/irst kind using orthogonality and integral transforms.}


\Classification{26A42; 33C05;  33C10; 33C45; 35A08}

\renewcommand{\thefootnote}{\arabic{footnote}}
\setcounter{footnote}{0}

\section{Introduction}
\label{Introduction}

In \cite{Cohl12pow} and~\cite{CTRS}
(see also~\cite{Cohlerratum12}),
we
present some def\/inite integral and inf\/inite series addition
theorems which arise from expanding fundamental solutions of
elliptic equations on $\R^d$ in axisymmetric coordinate systems
which separate Laplace's equation. We utilize orthogonality
and integral transforms to obtain new def\/inite integrals from some
of these addition theorems.

\section{Def\/inite integrals from integral transforms}
\subsection{Application of Hankel's transform}
We use the following result where for $x\in(0,\infty)$ we def\/ine
\[
F(r\pm 0):=\lim_{x\to r\pm}F(x);
\]
see 
\cite[p.~456]{Watson}:

\begin{theorem}
\label{2:Hankel}
Let $F:(0,\infty)\to\C$ be such that
\begin{gather}\label{2:cond}
\int_0^\infty \sqrt{x}\,|F(x)|\,dx<\infty,
\end{gather}
and let $\nu \ge -\frac12$. Then{\samepage
\begin{gather}\label{2:Hankel2}
 \frac12(F(r+0)+F(r-0))=\int_0^\infty uJ_\nu(ur)\int_0^\infty xF(x)J_\nu(ux)\,dx\,du
\end{gather}
provided that the positive number $r$ lies inside an interval in
which $F(x)$ has finite variation.}
\end{theorem}

As an illustration for the method of integral transforms, we give
the following example.  According to 
\cite[(13.22.2)]{Watson} (see also 
\cite[(6.612.3)]{Grad}), we have for $\operatorname{Re} a>0$,
$b,c>0$, $\operatorname{Re} \nu>-\frac12,$ then
\begin{gather}
\int_0^\infty e^{-ka}J_\nu(kb) J_\nu(kc) dk= \frac{1}{\pi\sqrt{bc}}
Q_{\nu-1/2}\left( \frac{a^2+b^2+{c}^2}{2bc}\right),
\label{QJ}
\end{gather}
where $J_\nu:\C\setminus(-\infty,0]\to\C$, for order $\nu\in\C$ is the Bessel function
of the f\/irst kind def\/ined in \cite[(10.2.2)]{NIST} and
$Q_\nu^\mu:\C\setminus(-\infty,1]\to\C$ for $\nu+\mu\notin -\N$
with degree $\nu$ and order $\mu$, is the
associated Legendre function of the second kind def\/ined in
\cite[(14.3.7), \S~14.21]{NIST}.  The Legendre function
of the second kind
$Q_\nu:\C\setminus(-\infty,1]\to\C$ for $\nu\notin -\N$ is def\/ined in
terms of the zero-order associated Legendre function of the second kind
$Q_\nu(z):=Q_\nu^0(z)$.
If we apply Theorem~\ref{2:Hankel} to the function
$F:(0,\infty)\to\C$ def\/ined by
\begin{gather}
F(k):=\frac{\pi\sqrt{c}}{k}e^{-ka}J_\nu(kc),
\label{FQJE}
\end{gather}
then condition~\eqref{2:cond} is satisf\/ied.
If we use~(\ref{FQJE}) in~\eqref{2:Hankel2} then we obtain the
following result.
If $\operatorname{Re} a>0$, $c>0$, $\operatorname{Re} \nu>-\frac12$, then
\[
\int_0^\infty
J_\nu(kb)\,
Q_{\nu-1/2}\left(\frac{a^2+b^2+c^2}{2bc}\right)\!
\sqrt{b}\,db=\frac{\pi\sqrt{c}}{k}e^{-ka}J_\nu(kc),
\]
which is actually given in
 \cite[(2.18.8.11)]{Prud90}.

Hardy 
 \cite[(33.16)]{Hardy08} derives an interesting extension
of (\ref{QJ}) (see also 
\cite[p.~389]{Watson} and
\cite[p.~17]{Askey75}).
We apply the Whipple formula \cite[(14.9.17), \S~14.21]{NIST}
to Hardy's extension to obtain
\begin{gather}
 \int_0^\infty k e^{-ka}J_\nu(kb)J_\nu(kc)dk\nonumber\\[0.10cm]
\qquad=\frac{-2a}{\pi\sqrt{bc}
\left(a^2+(b+c)^2\right)^{1/2}
\left(a^2+(b-c)^2\right)^{1/2}}
Q_{\nu-1/2}^1\left(\frac{a^2+b^2+c^2}{2bc}\right),
\label{Qnu1kint}
\end{gather}
for $\operatorname{Re} a>0$, $b,c>0$,
$\operatorname{Re} \nu>-1$.
It is mentioned in \cite[p.~389]{Watson} that~(\ref{Qnu1kint}) can be derived from~(\ref{QJ}) by dif\/ferentiation with respect to~$a$.
Using this integral and Theorem~\ref{2:Hankel}, we prove the following
theorem.

\begin{theorem}
Let $\operatorname{Re} a>0$, $c>0$,
$\operatorname{Re} \nu>-1$. Then
\begin{gather*}
\int_0^\infty
\frac{J_\nu(kb)}
{
\left(a^2+(b+c)^2\right)^{1/2}
\left(a^2+(b-c)^2\right)^{1/2}
}
Q_{\nu-1/2}^1\!\left(\frac{a^2+b^2+c^2}{2bc}\right)\sqrt{b}\,db
=-\frac{\pi\sqrt{c}}{2a}e^{-ka}J_\nu(kc).
\end{gather*}
\end{theorem}

\begin{proof} By applying Theorem~\ref{2:Hankel} to the function
$F:(0,\infty)\to\C$ def\/ined by
\[
F(k):=-\frac{\pi\sqrt{c}}{2a}e^{-ka}J_\nu(kc),
\]
using
(\ref{Qnu1kint}),
we obtain the desired result.
\end{proof}

Now, we give another example of how an integral
expansion for a fundamental solution of Laplace's equation on~$\R^3$
in parabolic coordinates can be used to prove a new def\/inite integral.

\begin{theorem}
Let $m\in\N_0$, $\lambda'\in(0,\infty)$, $\mu,\mu'\in(0,\infty)$,
$\mu\ne \mu'$, $k\in(0,\infty)$. Then
\[
\int_0^\infty Q_{m-1/2}(\chi) J_m(k\lambda)\sqrt{\lambda}\,d\lambda=
2\pi\sqrt{\lambda'\mu\mu'} J_m(k\lambda') I_m(k\mu_<) K_m(k\mu_>),
\]
where
\[
\chi=\frac{4\lambda^2\mu^2+4{\lambda'}^2{\mu'}^2+(\lambda^2-{\lambda'}^2+{\mu'}^2-\mu^2)^2}
{8\lambda\lambda'\mu\mu'}>1
\]
and $\mu_\lessgtr:={\min \atop \max}\{\mu,\mu'\}$.
\end{theorem}

\begin{proof}
We apply Theorem \ref{2:Hankel} to the function $F:(0,\infty)\to\C$ def\/ined by
\[
F(k):=2\pi\sqrt{\lambda'\mu\mu'} J_m(k\lambda') I_m(k\mu_<) K_m(k\mu_>),
\]
where $I_\nu:\C\setminus(-\infty,0]\to\C,$ for $\nu\in\C$, is the modif\/ied
Bessel function of the f\/irst kind def\/ined in \cite[(10.25.2)]{NIST}
and $K_\nu:\C\setminus(-\infty,0]\to\C$, for $\nu\in\C$, is the modif\/ied
Bessel function of the second kind def\/ined in~\cite[(10.27.4)]{NIST}.
Again, we see that condition~\eqref{2:cond} is satisf\/ied. We
obtain the desired result from 
\cite[(41)]{CTRS},
namely, for $\lambda\in(0,\infty)$,
\begin{gather*}
\int_0^\infty J_m(k\lambda) J_m(k\lambda') I_m(k\mu_<) K_m(k\mu_>)
kdk= \frac{Q_{m-1/2}(\chi)}{2\pi\sqrt{\lambda\lambda'\mu\mu'}}.\tag*{\qed}
\end{gather*}
  \renewcommand{\qed}{}
\end{proof}

\subsection{Application of Fourier cosine transform}

\begin{theorem}
Let
$a,b\in(0,\infty)$
with $b\le a$,
$k\in(0,\infty)$,
$\operatorname{Re} \nu>-\frac12$.
Then
\begin{gather}
\int_0^\infty Q_{\nu-1/2}\left( \frac{a^2+{b}^2+z^2}{2ab}\right)
\cos(kz)\,dz =\pi\sqrt{ab}\,I_\nu(kb)K_\nu(ka).
\label{Qmumhalfcoskz}
\end{gather}
\end{theorem}

\begin{proof}
According to 
\cite[(6.672.4)]{Grad}
we have the integral relation
\[
\int_0^\infty I_\nu(kb) K_\nu(ka) \cos(kz)
dk=\frac{1}{2\sqrt{ab}}
Q_{\nu-1/2}\left(\frac{a^2+{b}^2+z^2}{2ab}\right),
\]
where $a,b\in(0,\infty)$ with $b<a,$ $z>0$,
$\operatorname{Re} \nu>-\frac12$.
We obtain the desired result from Theorem~\ref{2:Hankel} with
$F:(0,\infty)\to\C$ def\/ined such that
\[
F(k):=
\pi
\sqrt{\frac{ab}{k}}\,I_\nu(kb)K_\nu(ka)
\]
and $\nu=-\frac12$.   Furthermore, if one makes the replacement
$z\mapsto z/(\sqrt{2}a)$, $k\mapsto\sqrt{2}ka$ in
\cite[(7.162.6)]{Grad}, namely
\[
\int_0^\infty Q_{\nu-1/2}(1+z^2)\cos(kz)dz=\frac{\pi}{\sqrt{2}}
I_{\nu}
\left(\frac{k}{\sqrt{2}}\right)
K_{\nu}
\left(\frac{k}{\sqrt{2}}\right),
\]
where $k\in(0,\infty)$ and $\operatorname{Re} \nu>-\frac12,$ then we see
that~(\ref{Qmumhalfcoskz}) holds also for any $a,b\in(0,\infty)$
with $a=b$.
\end{proof}

\section{Def\/inite integrals from orthogonality relations}
\label{Definiteintegralsfromorthogonalityrelations}

\subsection{Degree orthogonality for associated Legendre functions\\
with integer degree and order}
\label{DegreeorthogonalityforassociatedLegendrefunctions}

We take advantage of the degree orthogonality relation for the Ferrers
function of the f\/irst kind with integer degree and order, namely
(cf.\ 
\cite[(7.112.1)]{Grad})
\begin{gather}
\int_0^\pi {\mathrm P}_n^m(\cos\theta){\mathrm P}_{n'}^m(\cos\theta)\sin\theta d\theta
=\frac{2}{2n+1}\frac{(n+m)!}{(n-m)!}\delta_{n,n'},
\label{orthoglegendrePn}
\end{gather}
where $m,n,n'\in\N_0$, and $m\le n$, $m\le n'$.  We are using the
associated Legendre function of the f\/irst kind (on-the-cut),
${\mathrm P}_\nu^\mu:(-1,1)\to\C,$ for $\nu,\mu\in\C$, the Ferrers
function of the f\/irst kind, which is def\/ined in \cite[(14.3.1)]{NIST}.

The following estimates for the Ferrers function of the f\/irst kind with
integer degree and order will be useful.
If $\theta\in[0,\pi]$ and $m,n\in\N_0$ then \cite[\S~5.3, (19)]{Schafke63}
\begin{gather}\label{3:est1}
|{\mathrm P}_n^m(\cos\theta)|\le \frac{(m+n)!}{n!}
\end{gather}
and if $\theta\in(0,\pi)$ then \cite[p.~203]{MOS}
\begin{gather}
\label{3:est2}
|{\mathrm P}_n^m(\cos\theta)|< 2\frac{(n+m)!}{n!}
\left({\pi n}\right)^{-1/2}(\csc \theta)^{m+1/2}.
\end{gather}
If $\mu\in\C$, $\xi>0$ are f\/ixed and $0\le \nu\to+\infty$, we also have the
following asymptotic formulas for the associated Legendre functions
\begin{gather}
P_\nu^\mu(\cosh \xi) =
(2\pi\sinh\xi)^{-1/2}
\frac{\Gamma(\nu+\mu+1)}{\Gamma(\nu+\frac32)}
e^{(\nu+\frac12)\xi}\big(1+O\big(\nu^{-1}\big)\big)
\label{3:asyP}
,\\
Q_\nu^\mu(\cosh \xi) =
\left(\frac{\pi}{2\sinh\xi}\right)^{1/2}
\frac{\Gamma(\nu+\mu+1)}{\Gamma(\nu+\frac32)}
e^{-(\nu+\frac12)\xi+i\pi\mu} \big(1+O\big(\nu^{-1}\big)\big),
\label{3:asyQ}
\\
P_\nu^\mu(i\sinh \xi) =
(2\pi\cosh\xi)^{-1/2}
\frac{\Gamma(\nu+\mu+1)}{\Gamma(\nu+\frac32)}
e^{(\nu+\frac12)\xi+i\pi\nu/2}\big(1+O\big(\nu^{-1}\big)\big)
\label{3:asyPi}
,\\
Q_\nu^\mu(i\sinh \xi) =
\left(\frac{\pi}{2\cosh\xi}\right)^{1/2}
\frac{\Gamma(\nu+\mu+1)}{\Gamma(\nu+\frac32)}
e^{-(\nu+\frac12)\xi
-i\pi(\nu+1)/2
+i\pi\mu
} \big(1+O\big(\nu^{-1}\big)\big).
\label{3:asyQi}
\end{gather}
These asymptotic formulae follow from representations of Legendre
functions by Gauss hypergeometric functions; see \cite[(8.1.1), (8.10.4--5)]{Abra}.

\begin{theorem}
Let $n,m\in\N_0$, with
$n\ge m,$ $\nu\in\C\setminus\{2m,2m+2,2m+4,\ldots\}$, $r,r'\in(0,\infty)$,
$r\ne r'$, $\theta'\in(0,\pi)$. Then
\begin{gather}
 \int_0^\pi \big(\chi^2-1\big)^{(\nu+1)/4}Q_{m-1/2}^{-(\nu+1)/2}(\chi){\mathrm P}_n^m(\cos\theta)
(\sin\theta)^{(\nu+2)/2}d\theta\nonumber\\
\qquad{} =
\frac{i\sqrt{\pi}}{2^{(\nu+1)/2}(\sin\theta')^{\nu/2}}
\left(\frac{r_>^2-r_<^2}{rr'}\right)^{(\nu+2)/2}
Q_n^{-(\nu+2)/2}\left(\frac{r^2+{r'}^2}{2rr'}\right){\mathrm P}_n^m(\cos\theta'),
\label{bignuresult}
\end{gather}
where
\begin{gather}
\chi=\frac{r^2+{r'}^2-2rr'\cos\theta\cos\theta'}{2rr'\sin\theta\sin\theta'}
\label{chisphR3}
\end{gather}
and $r_{\lessgtr}:={\min \atop \max}\{r,r^\prime\}.$
\end{theorem}

\begin{proof}  We start with the following addition theorem for the
associated Legendre function of the second kind (see 
\cite{Cohl12pow}),
namely for $\theta\in(0,\pi)$,
\begin{gather}
\big(\chi^2-1\big)^{(\nu+1)/4}(\sin\theta)^{\nu/2}
Q_{m-1/2}^{-(\nu+1)/2}(\chi)
=\frac{i\sqrt{\pi}}{2^{(\nu+3)/2}}(\sin\theta^\prime)^{-\nu/2}
\left(\frac{r_>^2-r_<^2}{rr^\prime}\right)^{(\nu+2)/2}
\nonumber\\
\qquad{} \times\sum_{n=m}^\infty
(2n+1)\frac{(n-m)!}{(n+m)!}
Q_n^{-(\nu+2)/2}\left(\frac{r^2+{r^\prime}^2}{2rr^\prime}
\right) {\mathrm P}_n^m(\cos\theta) {\mathrm P}_n^m(\cos\theta^\prime),
\label{addtheorem3ba}
\end{gather}
where $\chi>1 $ is given by~(\ref{chisphR3}).
By~\eqref{3:est1} and~\eqref{3:asyQ} the inf\/inite series is uniformly convergent for $\theta\in[0,\pi]$.
Therefore, if we multiply both
sides of~(\ref{addtheorem3ba}) by $\sin\theta\, {\mathrm P}_{n'}^m(\cos\theta)$, where $n'\in\N_0$
and integrate over $\theta\in(0,\pi)$ we obtain~(\ref{bignuresult}).
\end{proof}

\begin{corollary}
Let $n,m\in\N_0$ with $n\ge m$, $r,r'\in(0,\infty)$, $r\ne r'$, $\theta'\in(0,\pi)$.
Then
\[
\int_0^\pi Q_{m-1/2}(\chi){\mathrm P}_n^m(\cos\theta)\sqrt{\sin\theta}d\theta=
\frac{2\pi
\sqrt{\sin\theta'}
}{2n+1}
{\mathrm P}_n^m(\cos\theta')\left(\frac{r_<}{r_>}\right)^{n+1/2},
\]
where $\chi>1$ is given by~\eqref{chisphR3}.
\end{corollary}

\begin{proof}
Substitute $\nu=-1$ in
(\ref{bignuresult}) and use 
\cite[(14.5.17)]{NIST}.
\end{proof}

\begin{theorem}\label{prolate}
Let $m,n\in\N_0$ with
$n\ge m$,
$\sigma,\sigma'\in(0,\infty)$,  $\theta'\in(0,\pi)$.
Then
\begin{gather}
 \int_0^\pi Q_{m-1/2}(\chi) {\mathrm P}_n^m(\cos\theta)
\sqrt{\sin\theta}\,d\theta
=2\pi (-1)^m\frac{(n-m)!}{(n+m)!}\nonumber\\
\qquad {} \times\sqrt{\sinh\sigma\sinh\sigma'\sin\theta'}\,
{\mathrm P}_n^m(\cos\theta') P_n^m(\cosh\sigma_<)Q_n^m(\cosh\sigma_>),
\label{prolateint}
\end{gather}
where
\begin{gather}
\chi=\frac{\cosh^2\sigma+\cosh^2\sigma'-\sin^2\theta-\sin^2\theta'
-2\cosh\sigma\cosh\sigma'\cos\theta\cos\theta'}
{2\sinh\sigma\sinh\sigma'\sin\theta\sin\theta'}
\label{chiprolate}
\end{gather}
and
$\sigma_{\lessgtr}:={\min \atop \max}\{\sigma,\sigma^\prime\}$.
\end{theorem}

\begin{proof} We start with the following addition theorem
for the associated Legendre function of the second kind
(see 
\cite[(37)]{CTRS}), namely
\begin{gather}
 Q_{m-1/2}(\chi)=\pi(-1)^m
\sqrt{\sinh\sigma\sinh\sigma'\sin\theta\sin\theta'}
\sum_{n=m}^\infty
(2n+1)
\left[\frac{(n-m)!}{(n+m)!}\right]^2
\nonumber\\
\hphantom{Q_{m-1/2}(\chi)=}{}
\times
{\mathrm P}_n^m(\cos\theta){\mathrm P}_n^m(\cos\theta')P_n^m(\cosh\sigma_<)Q_n^m(\cosh\sigma_>),
\label{addtheoremprolate}
\end{gather}
where $\chi\ge 1$ is def\/ined by (\ref{chiprolate}).
Note that $\chi=1$ only if $\sigma=\sigma'$ and $\theta=\theta'$, and in that
case $Q_{m-1/2}(\chi)$ has a logarithmic singularity.
If $\sigma\ne \sigma'$ then \eqref{3:est1}, \eqref{3:asyP}, \eqref{3:asyQ} show that the series in \eqref{addtheoremprolate}
is uniformly convergent for $\theta\in[0,\pi]$.
Therefore, if we multiply
both sides of (\ref{addtheoremprolate}) by
$\sqrt{\sin\theta}\,{\mathrm P}_{n'}^m(\cos\theta)$ and
integrate over $\theta\in[0,\pi]$, then
by~(\ref{orthoglegendrePn}) we have obtained (\ref{prolateint}).
If $\sigma=\sigma'$ then one may use \eqref{3:est2}, \eqref{3:asyP}, \eqref{3:asyQ} and the orthogonality relation \eqref{orthoglegendrePn} to show that the series in \eqref{addtheoremprolate}
(as a series of functions in the variable $\theta\in(0,\pi)$) converges
in $L^2(0,\pi)$. That is
\begin{gather}\label{L2}
\sum_{n=m}^\infty
\left\{\frac{(n-m)!}{(n+m)!}\phi_n(\theta')
P_n^m(\cosh \sigma_<) Q_n^m(\cosh \sigma_>)\right\}^2 <\infty,
\end{gather}
where $\phi_n\in L^2(0,\pi)$ forms an orthonormal basis and is def\/ined as
\[
\phi_n(\theta)
:=\sqrt{\sin\theta} \sqrt{\frac{2n+1}{2}
\frac{(n-m)!}{(n+m)!}}{\rm P}_n^m(\cos\theta),
\]
for $n=m,m+1,\dots$.
Then by the asymptotics of $P_n^m$ and $Q_n^m$
(cf.~(\ref{3:asyP}), (\ref{3:asyQ}))
\[
P_n^m(\cosh \sigma_<) Q_n^m(\cosh \sigma_>)=
O\left(n^{2m-1}\right)\quad\text{as $n\to\infty$}.
\]
Also from the estimate of ${\rm P}_n^m$ (\ref{3:est2}),
$\phi_n(\theta')=O\left(1\right).$
Therefore,
\[
\frac{(n-m)!}{(n+m)!} \phi_n(\theta')P_n^m(\cosh \sigma_<)
Q_n^m(\cosh \sigma_>)=O\big(n^{-1}\big),
\]
and this implies \eqref{L2} because
$\sum\limits_{n=1}^\infty \frac{1}{n^2} <\infty.$
Therefore, we again obtain (\ref{prolateint}).
\end{proof}

\begin{theorem}\label{oblate}
Let $m,n\in\N_0$ with $0\le m\le n$,
$\sigma,\sigma'\in(0,\infty)$, $\theta'\in[0,\pi].$
Then
\begin{gather}
 \int_0^\pi Q_{m-1/2}(\chi) {\mathrm P}_n^m(\cos\theta)\sqrt{\sin\theta}d\theta
=2\pi i(-1)^m\frac{(n-m)!}{(n+m)!}\nonumber\\
\qquad{} \times\sqrt{\cosh\sigma\cosh\sigma'\sin\theta'}
{\mathrm P}_n^m(\cos\theta') P_n^m(i\sinh\sigma_<)Q_n^m(i\sinh\sigma_>),
\label{oblateint}
\end{gather}
where
\begin{gather}
\chi=\frac{\sinh^2\sigma+\sinh^2\sigma'+\sin^2\theta+\sin^2\theta'
-2\sinh\sigma\sinh\sigma'\cos\theta\cos\theta'}
{2\cosh\sigma\cosh\sigma'\sin\theta\sin\theta'}.
\label{chioblate}
\end{gather}
\end{theorem}

\begin{proof} We start with oblate spheroidal coordinates on $\R^3$, namely
\[
x = a\cosh\sigma\sin\theta\cos\phi,\qquad
y = a\cosh\sigma\sin\theta\sin\phi,\qquad
z = a\sinh\sigma\cos\theta,
\]
where $a>0$, $\sigma\in [0,\infty)$, $\theta\in [0,\pi]$,
$\phi\in[0,2\pi)$.
The reciprocal distance between two points $\bfx,\bfxp\in\R^3$ expanded in terms
of the separable harmonics in this coordinate system is given in 
\cite[(41), p.~218]{MacRobert47}, namely
\begin{gather*}
 \frac{1}{\|\bfx-\bfxp\|}=\frac{i}{a}
\sum_{n=0}^\infty (2n+1)\sum_{m=0}^n(-1)^m\epsilon_m
\left[\frac{(n-m)!}{(n+m)!}\right]^2\cos(m(\phi-\phi')) \\
\hphantom{\frac{1}{\|\bfx-\bfxp\|}=}{}
\times {\mathrm P}_n^m(\cos\theta){\mathrm P}_n^m(\cos\theta')P_n^m(i\sinh\sigma_<)Q_n^m(i\sinh\sigma_>),
\end{gather*}
where $\epsilon_m:=2-\delta_{m,0}$ is the Neumann factor
\cite[p.~744]{MorseFesh} commonly occurring in Fourier cosine
series, with $\sigma'\in[0,\infty),$
$\theta'\in[0,\pi],$
$\phi'\in[0,2\pi)$. Note that the corresponding expression given in 
\cite[\S~5.2]{CTRS} is given incorrectly
(see 
\cite{Cohlerratum12}).
By reversing the order of summations in the above expression and
comparing with the Fourier cosine expansion for the reciprocal
distance between two points, namely
\[
\frac{1}{\|\bfx-\bfxp\|}=\frac{1}{\pi a
\sqrt{\cosh\sigma\cosh\sigma'\sin\theta\sin\theta'}}
\sum_{m=0}^\infty \epsilon_m \cos(m(\phi-\phi'))
Q_{m-1/2}(\chi),
\]
where $\chi>1$ is given by (\ref{chioblate}), we obtain the following
addition theorem for the associated Legendre function of the second kind
\begin{gather}
 Q_{m-1/2}(\chi)=i\pi(-1)^m\sqrt{\cosh\sigma\cosh\sigma'\sin\theta\sin\theta'}
\sum_{n=m}^\infty
(2n+1)
\left[\frac{(n-m)!}{(n+m)!}\right]^2
\nonumber\\
\hphantom{Q_{m-1/2}(\chi)=}{}
\times
{\mathrm P}_n^m(\cos\theta){\mathrm P}_n^m(\cos\theta')P_n^m(i\sinh\sigma_<)Q_n^m(i\sinh\sigma_>).
\label{addtheoremoblate}
\end{gather}
If we multiply both sides of (\ref{addtheoremoblate}) by
$\sqrt{\sin\theta}\,{\mathrm P}_{n'}^m(\cos\theta)$ and integrate over
$\theta\in[0,\pi]$, then by (\ref{orthoglegendrePn}) we have
obtained~(\ref{oblateint}).  We justify the interchange of integral and inf\/inite sum
as before by using the asymptotic formulas \eqref{3:asyPi}, \eqref{3:asyQi}.
\end{proof}

\begin{theorem}
Let $m,n\in\N_0$ with $0\le m\le n$,
$\sigma,\sigma'\in(0,\infty)$, $\theta'\in(0,\pi)$.
Then
\begin{gather}
 \int_0^\pi Q_{m-1/2}(\chi) {\mathrm P}_n^m(\cos\theta)
\sqrt{\sin\theta}d\theta
=\frac{2\pi\sqrt{\sin\theta'}}{2n+1} {\mathrm P}_n^m(\cos\theta') e^{-(n+1/2)(\sigma_>-\sigma_<)},
\label{bisphereint}
\end{gather}
where if we define $s=\cosh\sigma$, $s'=\cosh\sigma'$, $\tau=\cos\theta$,
$\tau^\prime=\cos\theta'$, then
\begin{gather}
\chi=\frac{
\sin^2\theta(s^\prime-\tau^\prime)^2+\sin^2\theta^\prime(s-\tau)^2
+\bigl[(s^\prime-\tau^\prime)\sinh\sigma-(s-\tau)\sinh\sigma^\prime\bigr]^2}
{2\sin\theta\sin\theta^\prime(s-\tau)(s^\prime-\tau^\prime)}.
\label{chibisphere}
\end{gather}
\end{theorem}

\begin{proof} We start with bispherical coordinates on $\R^3$, namely
\[
x = \frac{a\sin\theta\cos\phi}{\cosh\sigma-\cos\theta},\qquad
y = \frac{a\sin\theta\sin\phi}{\cosh\sigma-\cos\theta},\qquad
z = \frac{a\sinh\sigma}{\cosh\sigma-\cos\theta},
\]
where $a>0$, $\sigma\in [0,\infty)$, $\theta\in [0,\pi]$,
$\phi\in[0,2\pi)$.
The reciprocal distance between two points $\bfx,\bfxp\in\R^3$ expanded in terms
of the separable harmonics in this coordinate system is given in 
\cite[(9), p.~222]{MacRobert47}, namely
\begin{gather*}
 \frac{1}{\|\bfx-\bfxp\|}=\frac{1}{a}
\sqrt{(\cosh\sigma-\cos\theta)(\cosh\sigma'-\cos\theta')}
\sum_{n=0}^\infty e^{-(n+1/2)(\sigma_>-\sigma_<)} \\
\hphantom{\frac{1}{\|\bfx-\bfxp\|}=}{}
\times\sum_{m=0}^n\epsilon_m
\frac{(n-m)!}{(n+m)!}{\mathrm P}_n^m(\cos\theta){\mathrm P}_n^m(\cos\theta')
\cos(m(\phi-\phi')) ,
\end{gather*}
where $\sigma'\in[0,\infty)$,
$\theta'\in[0,\pi]$, $\phi'\in[0,2\pi)$.
By reversing the order of summations in the above expression and comparing
with the Fourier cosine expansion for the reciprocal distance between two
points, namely
\begin{gather*}
 \frac{1}{\|\bfx-\bfxp\|}=
\frac{\sqrt{(\cosh\sigma-\cos\theta)(\cosh\sigma'-\cos\theta')}}
{\pi a \sqrt{\sin\theta\sin\theta'}}
\sum_{m=0}^\infty \epsilon_m \cos(m(\phi-\phi'))
Q_{m-1/2}(\chi),
\end{gather*}
where $\chi\ge 1$
is given by (\ref{chibisphere}), we obtain the following
addition theorem for the associated Legendre function of the second kind
\begin{gather}
 Q_{m-1/2}(\chi)=\pi\sqrt{\sin\theta\sin\theta'}
\sum_{n=m}^\infty
\frac{(n-m)!}{(n+m)!}e^{-(n+1/2)(\sigma_>-\sigma_<)}{\mathrm P}_n^m(\cos\theta){\mathrm P}_n^m(\cos\theta').
\label{addtheorembisphere}
\end{gather}
If $\sigma=\sigma'$ and
$\theta=\theta'$, then $\chi=1,$ and $Q_{m-1/2}(\chi)$ has a logarithmic
singularity.
Note that the corresponding expression given in
\cite[\S~6.1, (45)]{CTRS} is given
incorrectly (see 
\cite{Cohlerratum12}). If we multiply both
sides of~(\ref{addtheorembisphere}) by
$\sqrt{\sin\theta}\,{\mathrm P}_{n'}^m(\cos\theta)$ and
integrate over $\theta\in[0,\pi]$, then
by~(\ref{orthoglegendrePn}) we have obtained~(\ref{bisphereint}).
We justify the interchange of integral and inf\/inite sum in the
same way as in the proof of Theorem~\ref{prolate}.
\end{proof}

\subsection{Order orthogonality for associated Legendre functions\\
with integer degree and order}
\label{OrderorthogonalityforassociatedLegendrefunctions}

In this subsection we take advantage of the order orthogonality relation for the
Ferrers function of the f\/irst kind with integer degree and order
(cf. 
\cite[(14.17.8)]{NIST})
\begin{gather}
\int_0^\pi {\mathrm P}_n^m(\cos\theta){\mathrm P}_n^{m'}(\cos\theta)\frac{1}{\sin\theta}d\theta
=\frac{1}{m}\frac{(n+m)!}{(n-m)!}\delta_{m,m'},
\label{orderPorthogonality}
\end{gather}
with $m\ge 1$.

\begin{theorem}\label{orderm}
Let $m\in\N$, $n\in\N_0$ with $1\le m\le n$,
$\theta'\in[0,\pi]$, $\phi,\phi'\in[0,2\pi)$. Then
\[
\int_0^\pi P_n(\cos\gamma){\mathrm P}_n^m(\cos\theta)\frac{1}{\sin\theta}d\theta=
\frac{2}{m}{\mathrm P}_n^m(\cos\theta')\cos(m(\phi-\phi')),
\]
where
\[
\cos\gamma=\cos\theta\cos\theta'+\sin\theta\sin\theta'\cos(\phi-\phi').
\]
\end{theorem}

\begin{proof}  We start with the addition theorem for spherical
harmonics (cf. 
\cite[(14.18.1)]{NIST}), namely
\begin{gather}
P_n(\cos\gamma)=\sum_{m=-n}^n\frac{(n-m)!}{(n+m)!}{\mathrm P}_n^m(\cos\theta){\mathrm P}_n^m(\cos\theta')
e^{im(\phi-\phi')},
\label{addtheoremsph}
\end{gather}
where $P_n:\C\to\C$, for $n\in\N_0$, is the Legendre polynomial which can be def\/ined
in terms of the terminating Gauss hypergeometric series (see for instance
\cite[Chapters~15,~18]{NIST}) as follows
\[
P_n(z):={}_2F_1\left(\!\begin{array}{c}-n,n+1\\ 1\end{array};\frac{1-z}{2}\right).
\]
We then take advantage of the order orthogonality relation for the Ferrers
functions of the f\/irst kind with integer degree and order.
If we multiply both sides of~(\ref{addtheoremsph}) by
$(\sin\theta)^{-1}\,{\mathrm P}_n^{m'}(\cos\theta)$ and integrate
over $\theta\in(0,\pi),$ by using~(\ref{orderPorthogonality}) we obtain
the desired result.
\end{proof}

Theorem~\ref{orderm}, originating from~(\ref{addtheoremsph}),
is the only example of a def\/inite integral that we could f\/ind
using the order orthogonality relation for the Ferrers functions of the f\/irst kind~(\ref{orderPorthogonality}).  Therefore we highly suspect that this result
is previously known, and include it mainly for completeness sake.  It would
however be very interesting to f\/ind another example using this orthogonality relation.

\subsection{Orthogonality for Chebyshev polynomials of the f\/irst kind}
\label{OrthogonalityfromChebyshevpolynomialsofthefirstkind}

Here we take advantage of orthogonality from Chebyshev polynomials
of the f\/irst kind
(cf. 
\cite[\S~18.3]{NIST})
\begin{gather}
\int_0^\pi T_m(\cos\theta) T_n(\cos\theta) d\theta=\frac{\pi}{\epsilon_n}
\delta_{m,n},
\label{orthogcheby1st}
\end{gather}
where   $T_n:\C\to\C$, for $n\in\N_0,$ is the Chebyshev polynomial
of the f\/irst kind which can  be def\/ined in terms of the terminating
Gauss hypergeometric series (see \cite[Chapter~18]{NIST})
\[
T_n(z)= {}_2F_1\left(\!\begin{array}{c}-n,n\\[1mm]\frac12\end{array};\frac{1-z}{2}\right).
\]
The Chebyshev polynomials of the f\/irst kind satisfy the identity
\cite[(18.5.1)]{NIST}
\begin{gather*}
T_n(\cos\theta)=\cos(n\theta).
\end{gather*}

\begin{theorem}
Let $m,n\in\Z$,
$\sigma,\sigma'\in(0,\infty)$.
Then
\begin{gather}
 \int_0^\pi Q_{m-1/2}(\chi) \cos(n\psi)d\psi
=\pi(-1)^m\sqrt{\sinh\sigma\sinh\sigma'}
\nonumber\\
\hphantom{\int_0^\pi Q_{m-1/2}(\chi) \cos(n\psi)d\psi=}{}
\times\frac{\Gamma\left(n-m+\frac12\right)}
{\Gamma\left(n+m+\frac12\right)}
P_{n-1/2}^m(\cosh\sigma_<)
Q_{n-1/2}^m(\cosh\sigma_>),
\label{toroidint}
\end{gather}
where
\begin{gather}
\chi=\coth\sigma\coth\sigma'-\operatorname{csch} \sigma \operatorname{csch} \sigma'\cos\psi.
\label{chitor}
\end{gather}
\end{theorem}

\begin{proof}  We start with toroidal coordinates on $\R^3$, namely
\[
x =  \frac{a\sinh\sigma\cos\phi}{\cosh\sigma-\cos\psi} , \qquad
y =  \frac{a\sinh\sigma\sin\phi}{\cosh\sigma-\cos\psi} ,\qquad
z = \frac{a\sin\psi}{\cosh\sigma-\cos\psi},
\]
where $a>0$, $\sigma\in (0,\infty)$, $\psi,\phi\in[0,2\pi)$.
The reciprocal distance between two points $\bfx,\bfxp\in\R^3$ is given
algebraically by
\begin{gather*}
 \frac{1}{\|\bfx-\bfxp\|}=\frac{1}{a}
\sqrt{
\frac{(\cosh\sigma-\cos\psi)(\cosh\sigma'-\cos\psi')}{2\sinh\sigma\sinh\sigma'}
}\\
\hphantom{\frac{1}{\|\bfx-\bfxp\|}=}{}
\times\left[
\frac{\cosh\sigma\cosh\sigma'-\cos(\psi-\psi')}{\sinh\sigma\sinh\sigma'}
-\cos(\phi-\phi')
\right]^{-1/2},
\end{gather*}
where $(\sigma',\psi',\phi')$ are the toroidal coordinates corresponding
to the point $\bfxp$.
Using Heine's recipro\-cal square root identity
(see for instance 
\cite[(3.11)]{CohlDominici})
\[
\frac{1}{\sqrt{z-x}}=\frac{\sqrt{2}}{\pi}
\sum_{m=0}^\infty \epsilon_m Q_{m-1/2}(z) T_m(x),
\]
where $z>1$ and $x\in[-1,1]$,
we can obtain a Fourier cosine series
representation for the reciprocal distance between two points in toroidal coordinates on $\R^3$, namely
\begin{gather*}
 \frac{1}{\|\bfx-\bfxp\|}=\frac{1}{\pi a}
\sqrt{\frac{(\cosh\sigma-\cos\psi)(\cosh\sigma'-\cos\psi')}
{\sinh\sigma\sinh\sigma'}}
\sum_{m=0}^\infty\epsilon_m\cos(m(\phi-\phi')) Q_{m-1/2}(\chi),
\end{gather*}
where $\chi>1$ is given by (\ref{chitor}).
We can further expand the
associated Legendre function of the second kind using the following addition theorem
(cf.~\cite[(8.795.2)]{Grad})
\begin{gather}
 Q_{m-1/2}(\chi)=(-1)^m\sqrt{\sinh\sigma\sinh\sigma'}
\sum_{n=0}^\infty\epsilon_n\cos(n(\psi-\psi'))
\nonumber\\
\hphantom{Q_{m-1/2}(\chi)=}{}
\times
\frac{\Gamma\left(n-m+\frac12\right)}{\Gamma\left(n+m+\frac12\right)}
P_{n-1/2}^m(\cosh\sigma_<) Q_{n-1/2}^m(\cosh\sigma_>).
\label{addtheoremtoroidal}
\end{gather}
Note that with the above addition theorem, we have the expansion of the
reciprocal distance between two points in terms of the separable harmonics in
toroidal coordinates
\begin{gather*}
 \frac{1}{\|\bfx-\bfxp\|}=\frac{1}{\pi a}
\sqrt{(\cosh\sigma-\cos\psi)(\cosh\sigma'-\cos\psi')}
\sum_{m=0}^\infty (-1)^m\epsilon_m\cos(m(\phi-\phi'))
\\
\hphantom{\frac{1}{\|\bfx-\bfxp\|}=}{}
\times
\sum_{n=0}^\infty\epsilon_n\cos(n(\psi-\psi'))
\frac{\Gamma\left(n-m+\frac12\right)}{\Gamma\left(n+m+\frac12\right)}
P_{n-1/2}^m(\cosh\sigma_<)Q_{n-1/2}^m(\cosh\sigma_>)
\end{gather*}
(see also 
 \cite[\S~6.2]{CTRS} and
 \cite{Cohlerratum12}).
If we relabel $\psi-\psi'\mapsto\psi$ and multiply both sides of
(\ref{addtheoremtoroidal}) by $\cos(n\psi)$ and integrate over
$\psi\in[0,\pi]$, then by~(\ref{orthogcheby1st}) we have obtained~(\ref{toroidint}). The interchange of inf\/inite
sum and integral is justif\/ied by~(\ref{3:asyP}),~(\ref{3:asyQ}).
\end{proof}

\subsection*{Acknowledgements}
This work was conducted while H.S.~Cohl was a National Research Council
Research Postdoctoral Associate in the Information Technology Laboratory at the
National Institute of Standards and Technology, Gaithersburg, Maryland, USA.
The authors would also like to acknowledge two anonymous referees
whose comments helped improve this paper.

\pdfbookmark[1]{References}{ref}
\LastPageEnding


\begin{thebibliography}{99}
\footnotesize\itemsep=0pt

\bibitem{Abra}
Abramowitz M., Stegun I.A., Handbook of mathematical functions with formulas,
  graphs, and mathematical tables, \textit{National Bureau of Standards Applied
  Mathematics Series}, Vol.~55, U.S. Government Printing Of\/f\/ice, Washington,
  D.C., 1964.

\bibitem{Askey75}
Askey R., Orthogonal polynomials and special functions, Society for Industrial
  and Applied Mathematics, Philadelphia, Pa., 1975.

\bibitem{Cohl12pow}
Cohl H.S., Fourier, Gegenbauer and Jacobi expansions for a power-law
  fundamental solution of the polyharmonic equation and polyspherical addition
  theorems, \href{http://arxiv.org/abs/1209.6047}{arXiv:1209.6047}.

\bibitem{Cohlerratum12}
Cohl H.S., Erratum: Developments in determining the gravitational potential
  using toroidal functions, \href{http://dx.doi.org/10.1002/asna.201211723}{\textit{Astronom. Nachr.}} \textbf{333} (2012),
  784--785.

\bibitem{CohlDominici}
Cohl H.S., Dominici D.E., Generalized {H}eine's identity for complex {F}ourier
  series of binomials, \href{http://dx.doi.org/10.1098/rspa.2010.0222}{\textit{Proc.~R. Soc. Lond. Ser.~A}} \textbf{467} (2011),
  333--345, \href{http://arxiv.org/abs/0912.0126}{arXiv:0912.0126}.

\bibitem{CTRS}
Cohl H.S., Tohline J.E., Rau A.R.P., Srivastava H.M., Developments in
  determining the gravitational potential using toroidal functions,
  \href{http://dx.doi.org/10.1002/1521-3994(200012)321:5/6<363::AID-ASNA363>3.0.CO;2-X}{\textit{Astronom. Nachr.}} \textbf{321} (2000), 363--372.

\bibitem{Grad}
Gradshteyn I.S., Ryzhik I.M., Table of integrals, series, and products, seventh
  ed., Elsevier/Academic Press, Amsterdam, 2007.

\bibitem{Hardy08}
Hardy G.H., Further researches in the theory of divergent series and integrals,
  \textit{Trans. Cambridge Philos. Soc.} \textbf{21} (1908), 1--48.

\bibitem{MacRobert47}
MacRobert T.M., Spherical harmonics. {A}n elementary treatise on harmonic
  functions with applications, 2nd ed., Methuen \& Co. Ltd., London, 1947.

\bibitem{MOS}
Magnus W., Oberhettinger F., Soni R.P., Formulas and theorems for the special
  functions of mathematical physics, 3rd ed., \textit{Die Grundlehren der
  mathematischen Wissenschaften}, Bd.~52, Springer-Verlag, New York, 1966.

\bibitem{MorseFesh}
Morse P.M., Feshbach H., Methods of theoretical physics, Vols.~1,~2,
  McGraw-Hill Book Co. Inc., New York, 1953.

\bibitem{NIST}
Olver F.W.J., Lozier D.W., Boisvert R.F., Clark C.W. (Editors), N{IST} handbook
  of mathematical functions, Cambridge University Press, Cambridge, 2010.

\bibitem{Prud90}
Prudnikov A.P., Brychkov Y.A., Marichev O.I., Integrals and series.
  {V}ol.~3.~More special functions, Gordon and Breach Science Publishers, New
  York, 1990.

\bibitem{Schafke63}
Sch{\"a}fke F.W., Einf\"uhrung in die {T}heorie der speziellen {F}unktionen der
  mathematischen {P}hysik, \textit{Die Grundlehren der mathematischen Wissenschaften},
  Bd.~118, Springer-Verlag, Berlin, 1963.

\bibitem{Watson}
Watson G.N., A treatise on the theory of {B}essel functions, 2nd ed., \textit{Cambridge
  Mathematical Library}, Cambridge University Press, Cambridge, 1944.

\end{thebibliography}
\end{document}